\documentclass[a4paper]{IEEEtran}

\usepackage{graphicx}
\usepackage[latin1]{inputenc}
\usepackage{epstopdf}
\usepackage[cmex10]{amsmath}
\usepackage{amsfonts}
\usepackage{amssymb}
\usepackage{mathrsfs}
\usepackage{verbatim}
\usepackage{array}
\usepackage{dcolumn}
\usepackage{color}
\usepackage{arydshln}
\usepackage{url}
\usepackage{mathrsfs}
\usepackage{algorithm, algorithmic}
\usepackage{relsize}
\usepackage{changepage}
\usepackage[bb=boondox]{mathalfa}
\usepackage[dvipsnames]{xcolor}
\usepackage{subfigure}
\usepackage{scalerel}
\usepackage{threeparttable}
\usepackage{tabularx}
\usepackage{multirow}
\usepackage{multicol}
\usepackage{makecell}

\usepackage{geometry}
\newcommand*{\pdot}{\mathbin{\scalerel*{\boldsymbol\odot}{\circ}}}

\allowdisplaybreaks
\usepackage{amssymb}
\usepackage{nomencl}
\usepackage{etoolbox}

\DeclareGraphicsExtensions{.eps}
\AppendGraphicsExtensions{.pdf}

\begin{document}

\title{Decision-Dependent Uncertainty-Aware Distribution System Planning Under Wildfire Risk}

\author{Felipe~Pianc\'o,~\IEEEmembership{Student~Member,~IEEE,} Alexandre Moreira, \IEEEmembership{Member,~IEEE,} Bruno Fanzeres, \IEEEmembership{Member,~IEEE,} Ruiwei Jiang, \IEEEmembership{Member,~IEEE,} Chaoyue Zhao, \IEEEmembership{Member,~IEEE,} and Miguel Heleno, \IEEEmembership{Senior Member,~IEEE}}
 
\maketitle

\begin{abstract}
The interaction between power systems and wildfires can be dangerous and costly. Damaged structures, load shedding, and high operational costs are potential consequences when the grid is unprepared. In fact, the operation of distribution grids can be liable for the outbreak of wildfires when extreme weather conditions arise. Within this context, investment planning should consider the impact of operational actions on the uncertainty related to wildfires that can directly affect line failure likelihood. Neglecting this can compromise the cost-benefit evaluation in planning system investments for wildfire risk. In this paper, we propose a decision-dependent uncertainty (DDU) aware methodology that provides the optimal portfolio of investments for distribution systems while considering that high power-flow levels through line segments in high-threat areas can ignite wildfires and, therefore, increase the probability of line failures. The methodology identifies the best combination of system upgrades (installation of new lines, hardening existing lines, and placement of switching devices) to provide the necessary leeway to operate the distribution system under wildfire-prone conditions. Our case study demonstrates that by modeling the DDU relationship between power flow prescriptions and line failures, investment decisions are more accurate and better prepare the grid infrastructure to deal with wildfire risk.

\end{abstract}

\begin{IEEEkeywords}
	Power Systems Investment Planning, Wildfire Resilience, Optimization Under Uncertainty, Decision-Dependent Uncertainty, Topology Reconfiguration
\end{IEEEkeywords}

\section*{Nomenclature}\label{Nomenclature}
\subsection*{Sets}
\begin{description} [\IEEEsetlabelwidth{100000}\IEEEusemathlabelsep]
    \item[$\mathcal{H}_{l}$] Set of hardening investments of a given line $l$.
    \item[$\mathcal{K}^{forb}$] Set of rules of topology forbidden patterns.
    \item[$\mathcal{L}$] Set of lines.
    \item[$\mathcal{L}^{e}$] Set of existing lines.
    \item[$\mathcal{L}^{c}$] Set of candidate lines.
    \item[$\mathcal{L}^{sw}$] Set of switchable line segments.
    \item[$\mathcal{L}^{sw,e}$] Set of existing switchable lines.
    \item[$\mathcal{L}^{sw,c}$] Set of candidates lines to become switchable.
    \item[$\mathcal{L}^{ha}$] Set of lines that have hardening investment options.
    \item[$\mathcal{L}^{forb}_{k}$] Set of lines in the forbidden pattern rule $k$ in $\mathcal{K}^{forb}$.
    \item[$\mathcal{N}$] Set of buses.
    \item[$\mathcal{N}^{sub}$] Set of buses with substation.
    \item[$\mathcal{R}$] Set of representative days.
    \item[$\mathcal{T}$] Set of periods.
\end{description}

\subsection*{Parameters}
\begin{description} [\IEEEsetlabelwidth{100000}\IEEEusemathlabelsep]
    \item[$\beta$] Sensitivity of line failure probability to the scheduled active power flow.
    \item[$\gamma$] Estimated upper bound for the nominal line failure probability.
    \item[$C^{p+}$] Cost of active power surplus.
    \item[$C^{p-}$] Cost of active power loss.
    \item[$C^{q+}$] Cost of reactive power surplus.
    \item[$C^{q-}$] Cost of reactive power loss.
    \item[$C^{sw}$] Cost of line switching action.
    \item[$C^{sw,inv}$] Investment cost of turning a line to be switchable.
    \item[$C^{line,inv}_l$] Investment cost of constructing a line.
    \item[$C^{ha,inv}$] Investment cost of line $l$ hardening action $h$ in $\mathcal{H}_{l}$.
    \item[$C^{tr}$] Cost of active power from main transmission grid at the substations (used in third-level).
    \item[$D^{p}$] Active power demand at the buses.
    \item[$\overline{F}$] Maximum power flow at the lines.
    \item[$K$] Number of simultaneous lines that can go off defined in the security criterion.
    \item[$M$] Sufficiently large number (Big M).
    \item[$\overline{P}$] Maximum active power injection at the buses.
    \item[$PF$] Power factor at the buses.
    \item[$\overline{Q}$] Maximum reactive power injection at the substations.
    \item[$\underline{Q}$] Minimum reactive power injection at the substations.
    \item[$R$] Resistance of the lines.
    \item[$S$] Auxiliary matrix of second-level problem.
    \item[$\underline{V}$] Voltage lower bound at the buses.
    \item[$\overline{V}$] Voltage upper bound at the buses.
    \item[$V^{ref}$] Voltage reference at the substations.
    \item[$X$] Reactance of the lines.
    \item[$w$] Weight parameter of first-level objective function.
    \item[$w^{ha}$] Weight of hardening investment decrease of $\beta$.
    \item[$z^{init}$] Initial switching status of the lines.
\end{description}

\subsection*{Decision Variables}
\begin{description} [\IEEEsetlabelwidth{100000}\IEEEusemathlabelsep]
    \item[${\Delta}D^{p+}$] Amount of buses active power surplus.
    \item[${\Delta}D^{p-}$] Amount of buses active power loss.
    \item[${\Delta}D^{q+}$] Amount of buses reactive power surplus.
    \item[${\Delta}D^{q-}$] Amount of buses reactive power loss.
    \item[$\varphi$] Dual decision variable of the worst expected value of the third-level problem.
    \item[$\psi$] Dual decision variable of the worst expected value of third-level problem for every line.
    \item[$a$] Binary variable associated with line availability (1 if available, 0 otherwise).
    \item[$f^{p}$] Active power flow at the lines.
    \item[$f^{q}$] Reactive power flow at the lines.
    \item[$p$] Amount of active power injected at the substations.
    \item[$q$] Amount of reactive power injected at the substations.
    \item[$v^{\dagger}$] Squared voltage at the buses.
    \item[$y^{sw,inv}$] Binary decision variable indicating an investment action for a line to become switchable (1 if an investment was made, 0 otherwise).
    \item[$y^{line,inv}$] Binary decision variable indicating an investment action to construct a line (1 if constructed, 0 otherwise).
    \item[$y^{ha,inv}$] Binary decision variable indicating a line hardening investment action (1 if invested, 0 otherwise).
    \item[$z^{sw}$] Binary decision variable indicating a line switching action (1 if switched, 0 otherwise).
    \item[$z^{topo}$] Binary decision variable of line topology (1 if on, 0 otherwise).
\end{description}

\newpage
\section{Introduction} \label{sec:Introduction}

% ====== Wildfire potential and  a hint of the paper's objective
\IEEEPARstart{T}{he} potential for wildfires has increased worldwide and the impact of human activities on wild lands combined with extreme weather conditions constitute one of the major causes for this issue \cite{Jahn2022}. In the context of power systems, their interaction with wildfires can be particularly harmful given the wide presence of electrical networks and the high dependency on electricity in modern society. Apart from damaging different structures, wildfires can directly affect the operation of power systems. For instance, the sudden growth of temperature near the lines can reduce power lines' capacity temporarily, and the combination of heat and smoke can change the properties of the air gap, increasing the probability of short circuits \cite{Serrano2022}. On the other hand, whilst wildfires have the potential to impact power systems, the operation of these systems, combined with certain environmental conditions, can also start wildfires. Within this context, it is critical to plan strategic investments that can increase the flexibility of power grids to avoid igniting wildfires while properly considering the role of electrical networks in initiating these events.

\subsection{Problem Context and Literature Review} \label{subsec:LiteratureReview}

\begin{figure}[t]
\centerline{
    \subfigure[]{\includegraphics[width=0.24\textwidth]{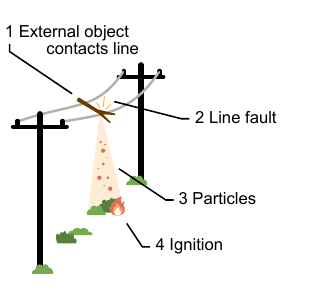}
    \label{subfig:PowerSystemWildfires_a}}
    \subfigure[]{\includegraphics[width=0.24\textwidth]{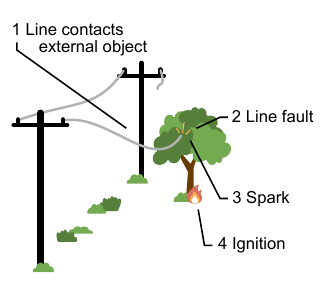}
    \label{subfig:PowerSystemWildfires_b}}}
\caption{Two examples of a wildfire started by power systems. Adapted from \cite{Jahn2022}.}
\label{fig:PowerSystemWildfires}
\end{figure}

% ======Wildfire and power system interaction
Wildfires induced by power systems can be started mainly in two manners: when external objects (e.g., tree parts) come into contact with the power line (Fig. \ref{subfig:PowerSystemWildfires_a}); or when the line itself comes into contact with nearby vegetation, equipment, or other lines (Fig. \ref{subfig:PowerSystemWildfires_b}) \cite{Wang2023}. In either case, the sequence of two events is necessary: a line fault followed by the ignition of nearby vegetation \cite{Muhs2020}. The ignition process occurs when the power flowing in the conductor finds an alternative path to the ground, closing a circuit. This action can directly start a fire when close to vegetation by creating an arc (Fig. \ref{subfig:PowerSystemWildfires_b}), or indirectly through molten metal particles, burning embers, or burning fluids (Fig. \ref{subfig:PowerSystemWildfires_a}). In this context, the higher the power flow, the higher the thermal stress of the lines, and the higher the chance of electric arcs to occur \cite{Muhs2020}. The duration and intensity of these arcs are directly related to ignition probability \cite{Victorian2011}.

% ====== Real cases of this interaction
Although power grid-induced ignitions are not the leading cause of wildfires, they are particularly threatening as they often occur in tandem with severe weather conditions, including high wind, low humidity, and high temperature \cite{Wang2023}. As extreme weather can comprise vast areas, it can lead to multiple line faults and multiple ignitions. In addition, certain weather conditions can create a mechanism of rapid-fire growth and ember showers, which might result in multiple large fires and consequent catastrophic losses \cite{Mitchell2013}. Some examples of wildfires initiated by power lines are (i) the ``Black Saturday'', which destroyed more than 400,000 ha and caused 173 human deaths in Victoria, Australia in 2009 \cite{Teague2010}; (ii) the ``Bastrop fire'', that resulted in the death of 2 people and the loss of hundreds of homes, in Texas, USA in 2011 \cite{Ridenour2012}; (iii) the ``Camp Fire'', which led to the death of 85 people in 2018, the deadliest wildfire in the history of California, USA \cite{Maranghides2021}; (iv) the ``Dixie fire'', which burned an area larger than 300,000 ha in California, USA in 2021 \cite{Calfire2019}. Indeed, in California, out of the 20 most destructive wildfires observed in the past decades, at least 9 were caused by the interaction between electrical infrastructure and environment.% \cite{Calfire2018}.

% ====== Resilience actions taken (California case)
Concerned about the numerous adverse fire incidents over the past two decades, electric utilities have started specific programs to tackle this issue in California. For instance, the San Diego Gas \& Electric (SDG\&E) has created numerous initiatives in this matter \cite{Udren2022}, including situational awareness, forecasting, and wildfire risk modeling programs. Furthermore, infrastructure improvements have also been implemented, encompassing conductors covering, equipment modernization, protection enhancement, and microgrid development \cite{Udren2022}. Southern California Edison (SCE) also developed special operational protocols and implemented vegetation clearance and conductor hardening. In addition, artificial intelligence and machine learning algorithms have helped the company to identify equipment defects via aerial inspections \cite{Chiu2022}. Similar actions have been taken by the Pacific Gas and Electric (PG\&E) \cite{PGE2023}. The company is one of the utilities in California that use the Public Safety Power Shutoff (PSPS) program, where, under specific weather conditions, power lines are actively de-energized in high-threat areas to prevent wildfire disruption. As a consequence, this program can cause power outages for customers during long periods \cite{Huang2023}. A more recent initiative from PG\&E, the so-called Enhanced Powerline Safety Settings (EPSS) program, is helping to reduce the frequency of PSPS events by deploying protection devices that can rapidly de-energize power lines once in contact with external objects, therefore decreasing the amount of energy released by a potential arc and thus reducing potential sources of wildfire ignition \cite{PGE2023}. Moreover, PG\&E also considers a major line undergrounding program in areas that are vulnerable to wildfires \cite{PGE2023}.

% ====== Operation and planning papers
Several technical and academic works have also provided insights into how to deal with the interaction between power systems and extreme events, including wildfire. From an operational perspective, numerous papers proposed different methodologies to support decision-making under these circumstances. For instance, \cite{Abdelmalak2021} proposes a probabilistic re-dispatch strategy for transmission systems using a Markov decision process that models state transitions while considering uncertainties in failure, load, and wildfires' spatial-temporal progression. In \cite{Huang2017}, the authors develop an integrated resilience response framework that comprises generation re-dispatch, transmission lines switching, and load shedding as preventive and emergency measures to deal with extreme events. The stochastic model proposed in \cite{Trakas2018} incorporates an exogenous wildfire progression to optimize emergency actions (including the dispatch of microturbines) to minimize load shedding in distribution grids. The work developed in \cite{Mehdi2023} proposes a risk-averse Information-Gap Decision (IGPD) method to support system operation while taking into consideration the impact of wildfires on the climate conditions in the immediate surroundings of distribution power lines, which can result in ampacity derates. In \cite{Lei2020} and \cite{Bie2017}, network configuration and microgrid formation are leveraged to alleviate the impact of extreme events in distribution systems. From a planning perspective, different approaches to optimize investments for power systems under extreme event threats can also be found in the literature. For example, in \cite{Bertoletti2022}, a transmission expansion problem is studied considering different fire-threat zones. Wind speed, vegetation, and humidity are the wildfire-related risk factors considered in \cite{Bayani2023} when selecting a portfolio of investments (including new lines and hardening actions) to enhance the resilience of the transmission system against wildfires. In \cite{Nagarajan2016}, a two-stage model is proposed to identify upgrade investments that can improve the capability of transmission systems to circumvent predefined adverse extreme scenarios. Furthermore, in \cite{Yanling2018} the authors designed a three-level optimization model to coordinate line hardening solutions and operational measures to protect the distribution grid against natural disasters in general. With a similar objective, in \cite{Byeon2020}, the authors propose a two-stage stochastic framework that models communication networks and optimizes different measures to harden the distribution system.

% ====== DDU papers
Despite the relevant contributions of the aforementioned approaches proposed by industry and academia, none of them has addressed the endogenous nature of wildfire uncertainty and operational decisions while planning distribution grid investments. Methodologically, the Decision-Dependent Uncertainty (DDU) framework models the endogenous nature of uncertainty realizations when the uncertainty dynamics are influenced by the solution assigned to all or to part of the decision variables of the problem \cite{Luo2020}. In the past years, a few scientific and technical works have considered the DDU framework in power-systems-related problems. In \cite{Shanshan2018}, for instance, the authors propose a hardening planning model where the uncertainty related to the impact of wind-induced adverse events depends on the line hardening decisions. Also, in \cite{Wenqian2023}, the authors propose a decision-dependent stochastic model for the operation and maintenance of transmission lines considering the occurrence of sandstorms and consequent line failures. In addition, in \cite{Zhang2023}, a two-stage stochastic program is developed to plan transmission hardening against typhoons. In this case, decision-dependency connects the relationship between line enhancement measures and the probability of line failure. Moreover, similarly to \cite{Zhang2023}, but in the context of earthquakes, the two-stage planning model proposed in \cite{Shi2023} considers a link between line failure uncertainty and hardening decisions. In \cite{Moreira2023}, we designed a methodology for tackling the operation of distribution grids capable of endogenously taking into account the impact of the line power flow in fault probabilities with a DDU framework. Nonetheless, to the best of our knowledge, there is still a need for a DDU-aware methodology that can plan investments to improve distribution system operation under wildfire risk.

\subsection{Contributions} \label{subsec:Contributions}

In this paper, we propose a methodology to devise expansion and investment plans for distribution systems while considering, in a DDU-aware framework, the double role of power lines as assets that can affect and be affected by the uncertainty realization related to wildfires. It is worth mentioning that, unlike the previously described works that account for DDU, our approach links investments and operational decisions with the very initiation of the extreme event as we consider that power lines in high-threat areas can ignite wildfires. Within this context, we consider the installation of new lines and the hardening of existing ones, as well as the strategic allocation of switches in the distribution grid, as candidate investments to be chosen by the system planner to increase the operational leeway under wildfire-prone conditions. Structurally, the proposed methodology in this paper falls into the class of a two-stage, three-level, Distributionally Robust Optimization (DRO) problem with DDU. We leverage the DRO framework as it is non-trivial to infer precise line failure probabilities. Structurally, in the first stage, the model determines the investment plan and the network topology configuration. In the second stage, the post-contingency multiperiod operation is addressed, where the probabilities are adjusted according to the first-stage information and the investments made. To summarize, the contributions of this paper are: 

\begin{enumerate}
    \item To design a distribution grid expansion planning methodology that accounts for the double-sided role of operational decisions in wildfire-prone seasons. More specifically, we formulate a DDU-aware two-stage DRO model that takes into account the impact of investments and power-flow levels on the wildfire uncertainty realization and the induced ambiguity state in the line failure probabilistic description to optimize investments aiming at providing distribution operators the necessary leeway to manage the system under wildfire-prone conditions.
    
    \item To propose a framework that explicitly maps the relationship between first-stage operational decisions and second-stage wildfire uncertainty impact to properly balance the cost/benefit trade-off when devising investment plans for distribution grids. For this purpose, we improve the standard adjustable DRO problem class to account for an ambiguity set of line-failure distributions that depends on investments and operational points. Structurally, the resulting mathematical formulation has a three-level system of nested optimization models.

    \item To construct a tailored iterative solution method capable of solving the proposed non-convex optimization problem in a finite number of steps. More specifically, we devise a set of reformulation procedures leveraging duality theory to develop an outer-approximation-based decomposition algorithm composed of tractable master- and sub-problems. Numerical experiments conducted in this work indicate that the proposed solution method can handle reasonable-sized instances in a practical computational time.
\end{enumerate}

\section{Problem Formulation} \label{sec:Problem}

\begin{figure}[t]
    \centering
    \includegraphics[width=0.48\textwidth,height=.5\textheight,keepaspectratio]{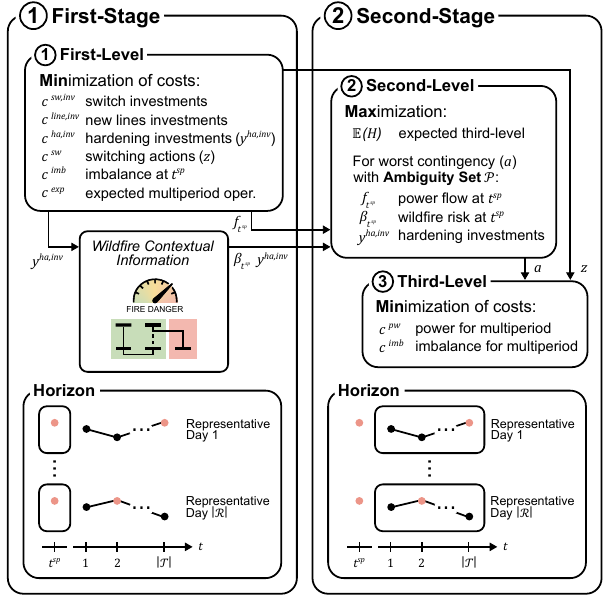}
    \caption{Proposed two-stage investment planning methodology as a three-level system of nested optimization models.}
    \label{fig:ModelDiagram}
\end{figure}

The main objective of the methodology proposed in this paper is to determine the optimal portfolio of investments for distribution grids while explicitly considering the impact of power-flow levels on the potential ignition of wildfires. Methodologically, as illustrated in Fig. \ref{fig:ModelDiagram}, we formulate this investment planning problem as a two-stage distributionally robust optimization model that accounts for decision-dependent uncertainty. The first stage defines investments and topology configuration to minimize overall system expenditures, including the worst-case expected second-stage costs under a decision-dependent ambiguity set. The second stage accounts for a multiperiod operation of the system. In the next subsections, we discuss in detail the mathematical formulation proposed in this paper.

\subsection{Investment Planning Model under Wildfire-Related DDU} \label{subsec:FirstLevel}

The investment planning model for distribution grids proposed in this work is presented in \eqref{1st_level_00}--\eqref{1st_level_33}.
\begin{align}
    & \underset{{\substack{
                           p_{b,t^{sp},r}, q_{b,t^{sp},r}, v^{\dagger}_{b,t^{sp},r}, f^{p}_{l,t^{sp},r}, f^{q}_{l,t^{sp},r}, \\
                           {\Delta}D^{p+}_{b,t^{sp},r}, {\Delta}D^{p-}_{b,t^{sp},r}, {\Delta}D^{q+}_{b,t^{sp},r}, {\Delta}D^{q-}_{b,t^{sp},r}, \\ 
                           z^{topo}_{l,r}, z^{sw}_{l,r}, y^{sw,inv}_{l}, y^{line,inv}_{l}, y^{ha,inv}_{l,h}
                           }}}{\text{Minimize}} 
      \hspace{0.00cm} \sum_{l \in {\mathcal{L}}^{sw,c}} C^{sw,inv}_{l} y^{sw,inv}_{l} \notag \\
    & \hspace{0.00cm} + \sum_{l \in {\mathcal{L}}^{c}} C^{line,inv}_{l} y^{line,inv}_{l} + \sum_{l \in \mathcal{L}^{ha}} \sum_{h \in \mathcal{H}_{l}} C^{ha,inv}_{l,h} y^{ha,inv}_{l,h} \notag \\
    & \hspace{0.00cm} + \sum_{r \in \mathcal{R}} w_{r} \Bigg( \sum_{l \in {\mathcal{L}}^{sw}} C^{sw}_{l} z^{sw}_{l,r} + \sum_{b \in \mathcal{N}} \Big( C^{p+} {\Delta}D^{p+}_{b,t^{sp},r} +  \notag \\
    & \hspace{1.70cm} C^{p-} {\Delta}D^{p-}_{b,t^{sp},r} + C^{q+} {\Delta}D^{q+}_{b,t^{sp},r} + C^{q-} {\Delta}D^{q-}_{b,t^{sp},r} \Big) \notag \\
    & \hspace{1.60cm} + \sup_{\mathcal{Q} \in \mathcal{P}_{r}(\boldsymbol{f}_{r}^{p},\boldsymbol{y}^{ha})} \mathbb{E}_\mathcal{Q} \big[ H_{r}(\boldsymbol{z}^{topo}_{r},\boldsymbol{a}_{r}) \big] \Bigg) \label{1st_level_00} \\ \notag
    & \text{subject to:} \notag \\ 
    & \hspace{0.00cm} p_{b,t^{sp},r} + \sum_{l \in \mathcal{L}|to(l)=b} f^{p}_{l,t^{sp},r} - \sum_{l \in \mathcal{L}|fr(l)=b} f^{p}_{l,t^{sp},r} - D^{p}_{b,t^{sp},r} \notag \\
    & \hspace{1.0cm} - {\Delta}D^{p+}_{b,t^{sp},r} + {\Delta}D^{p-}_{b,t^{sp},r} = 0 ; \forall b \in \mathcal{N}^{sub}, r \in \mathcal{R} \label{1st_level_01} \\
    & \hspace{0.00cm} q_{b,t^{sp},r} + \sum_{l \in \mathcal{L}|to(l)=b} f^{q}_{l,t^{sp},r} - \sum_{l \in \mathcal{L}|fr(l)=b} f^{q}_{l,t^{sp},r} \notag \\
    & \hspace{0.10cm} - \tan(\arccos(PF_{b})) D^{p}_{b,t^{sp},r} - {\Delta}D^{q+}_{b,t^{sp},r} + {\Delta}D^{q-}_{b,t^{sp},r}  \notag \\
    & \hspace{4.65cm} = 0; \forall b \in \mathcal{N}^{sub}, r \in \mathcal{R} \label{1st_level_02} \\
    & \hspace{0.00cm} \sum_{l \in \mathcal{L}|to(l)=b} f^{p}_{l,t^{sp},r} - \sum_{l \in \mathcal{L}|fr(l)=b} f^{p}_{l,t^{sp},r} - D^{p}_{b,t^{sp},r} \notag \\
    & \hspace{0.35cm} - {\Delta}D^{p+}_{b,t^{sp},r} + {\Delta}D^{p-}_{b,t^{sp},r} = 0 ; \forall b \in \mathcal{N} \setminus \mathcal{N}^{sub}, r \in \mathcal{R} \label{1st_level_03} \\
    & \hspace{0.00cm} \sum_{l \in \mathcal{L}|to(l)=b} f^{q}_{l,t^{sp},r} - \sum_{l \in \mathcal{L}|fr(l)=b} f^{q}_{l,t^{sp},r} \notag \\
    & \hspace{0.10cm} - \tan(\arccos(PF_{b})) D^{p}_{b,t^{sp},r} - {\Delta}D^{q+}_{b,t^{sp},r} + {\Delta}D^{q-}_{b,t^{sp},r} \notag \\
    & \hspace{4.00cm}  = 0 ; \forall b \in \mathcal{N} \setminus \mathcal{N}^{sub}, r \in \mathcal{R} \label{1st_level_04} \\
    & \hspace{0.00cm} v^{\dagger}_{to(l),t^{sp},r} - v^{\dagger}_{fr(l),t^{sp},r} + 2(R_{l} f^{p}_{l,t^{sp},r} + X_{l} f^{q}_{l,t^{sp},r}) \notag \\
    & \hspace{3.45cm} \leq (1 - z^{topo}_{l,r})M ; \forall l \in \mathcal{L}, r \in \mathcal{R} \label{1st_level_05} \\
    & \hspace{0.00cm} v^{\dagger}_{fr(l),t^{sp},r} - v^{\dagger}_{to(l),t^{sp},r} - 2(R_{l} f^{p}_{l,t^{sp},r} + X_{l} f^{q}_{l,t^{sp},r}) \notag \\
    & \hspace{3.45cm} \leq (1 - z^{topo}_{l,r})M ; \forall l \in \mathcal{L}, r \in \mathcal{R} \label{1st_level_06} \\
    & \hspace{0.00cm} f^{q}_{l,t^{sp},r} - \cot \Big( \Big( \frac{1}{2} - e  \Big) \frac{\pi}{4} \Big) \Big( f^{p}_{l,t^{sp},r} - \cos \Big( e \frac{\pi}{4} \Big) \overline{F}_l \Big) \notag \\
    & \hspace{0.8cm} - \sin \Big( e \frac{\pi}{4} \Big) \overline{F}_l \leq 0 ; \forall l \in \mathcal{L}, e \in \{1,\dots,4\}, r \in \mathcal{R} \label{1st_level_07} \\
    & \hspace{0.00cm} - f^{q}_{l,t^{sp},r} - \cot \Big( \Big( \frac{1}{2} - e  \Big) \frac{\pi}{4} \Big) \Big( f^{p}_{l,t^{sp},r} - \cos \Big( e \frac{\pi}{4} \Big) \overline{F}_l \Big) \notag \\
    & \hspace{0.8cm} - \sin \Big( e \frac{\pi}{4} \Big) \overline{F}_l \leq 0 ; \forall l \in \mathcal{L}, e \in \{1,\dots,4\}, r \in \mathcal{R} \label{1st_level_08} 
    \\  
    & \hspace{0.00cm} - z^{topo}_{l,r} \overline{F}_{l} \leq f^{p}_{l,t^{sp},r} \leq z^{topo}_{l,r} \overline{F}_{l} ; \forall l \in \mathcal{L}, r \in \mathcal{R} \label{1st_level_09} \\ 
    & \hspace{0.00cm} - z^{topo}_{l,r} \overline{F}_{l} \leq f^{q}_{l,t^{sp},r} \leq z^{topo}_{l,r} \overline{F}_{l} ; \forall l \in \mathcal{L}, r \in \mathcal{R} \label{1st_level_10} 
    \\
    & \hspace{0.00cm} \underline{V}^{2}_{b} \leq v^{\dagger}_{b,t^{sp},r} \leq \overline{V}^{2}_{b} ; \forall b \in \mathcal{N} \setminus \mathcal{N}^{sub}, r \in \mathcal{R} \label{1st_level_11} \\
    & \hspace{0.00cm} v^{\dagger}_{b,t^{sp},r} = V^{ref^2} ; \forall b \in \mathcal{N}^{sub}, r \in \mathcal{R} \label{1st_level_12}  \\
    & \hspace{0.00cm} 0 \leq p_{b,t^{sp},r} \leq \overline{P}_{b} ; \forall b \in \mathcal{N}^{sub}, r \in \mathcal{R} \label{1st_level_13} \\
    & \hspace{0.00cm} \underline{Q}_{b} \leq q_{b,t^{sp},r} \leq \overline{Q}_{b} ; \forall b \in \mathcal{N}^{sub}, r \in \mathcal{R} \label{1st_level_14} \\
    & \hspace{0.00cm} 0 \leq {\Delta}D^{p+}_{b,t^{sp},r} \leq D^{p}_{b,t^{sp},r} ; \forall b \in \mathcal{N}, r \in \mathcal{R} \label{1st_level_15} \\
    & \hspace{0.00cm} 0 \leq {\Delta}D^{p-}_{b,t^{sp},r} \leq D^{p}_{b,t^{sp},r} ; \forall  b \in \mathcal{N}, r \in \mathcal{R} \label{1st_level_16} \\
    & \hspace{0.00cm} 0 \leq {\Delta}D^{q+}_{b,t^{sp},r} \leq \tan(\arccos(PF_{b})) D^{p}_{b,t^{sp},r} ; \notag \\
    & \hspace{5.70cm} \forall b \in \mathcal{N}, r \in \mathcal{R} \label{1st_level_17} \\
    & \hspace{0.00cm} 0 \leq {\Delta}D^{q-}_{b,t^{sp},r} \leq \tan(\arccos(PF_{b})) D^{p}_{b,t^{sp},r} ; \notag \\
    & \hspace{5.70cm} \forall b \in \mathcal{N}, r \in \mathcal{R} \label{1st_level_18} \\
    & \hspace{0.00cm} z^{sw}_{l,r} = 0 ; \forall l \in \mathcal{L} \setminus \mathcal{L}^{sw}, r \in \mathcal{R} \label{1st_level_19} \\ 
    & \hspace{0.00cm} y^{sw,inv}_{l} = 0 ; \forall l \in (\mathcal{L}^{sw,e}) \cup (\mathcal{L} \setminus \mathcal{L}^{sw}) \label{1st_level_20} \\ 
    & \hspace{0.00cm} y^{line,inv}_{l} = 0 ; \forall l \in \mathcal{L}^{e} \label{1st_level_21} \\  
    & \hspace{0.00cm} z^{topo}_{l,r} = 1 ; \forall l \in \mathcal{L}^{e} \setminus \mathcal{L}^{sw}, r \in \mathcal{R} \label{1st_level_22} \\ 
    & \hspace{0.00cm} z^{sw}_{l,r} \leq y^{sw,inv}_{l} ; \forall l \in \mathcal{L}^{sw,c}, r \in \mathcal{R} \label{1st_level_23} \\ 
    & \hspace{0.00cm} y^{sw,inv}_{l} \leq y^{line,inv}_{l} ; \forall l \in (\mathcal{L}^{c}) \cap (\mathcal{L}^{sw,c}) \label{1st_level_24} \\     
    & \hspace{0.00cm} z^{sw}_{l,r} \geq + z^{topo}_{l,r} - z^{init}_{l,r} ; \forall l \in \mathcal{L}^e, r \in \mathcal{R} \label{1st_level_25} \\ 
    & \hspace{0.00cm} z^{sw}_{l,r} \geq - z^{topo}_{l,r} + z^{init}_{l,r} ; \forall l \in \mathcal{L}^e, r \in \mathcal{R} \label{1st_level_26} \\
    & \hspace{0.00cm} y^{line,inv}_{l} - y^{sw,inv}_{l} \leq z^{topo}_{l,r} ; \forall l \in \mathcal{L}^{c} \setminus \mathcal{L}^{sw,e}, r \in \mathcal{R} \label{1st_level_27} \\
    & \hspace{0.00cm} z^{topo}_{l,r} \leq y^{line,inv}_{l} ; \forall l \in \mathcal{L}^{c}, r \in \mathcal{R} \label{1st_level_28} \\ 
    & \hspace{0.00cm} \sum_{l \in \mathcal{L}^{forb}_{k}} z^{topo}_{l,r} \leq |\mathcal{L}^{forb}_{k}| - 1 ; \forall k \in \mathcal{K}^{forb}, r \in \mathcal{R} \label{1st_level_29} \\
    & \hspace{0.00cm} \sum_{h \in \mathcal{H}_{l}} y^{ha,inv}_{l,h} \leq 1 ; \forall l \in \mathcal{L}^{ha} \label{1st_level_30} \\
    & \hspace{0.00cm} z^{topo}_{l,r}, z^{sw}_{l,r} \in \{0,1\} ; \forall l \in \mathcal{L}, r \in \mathcal{R} \label{1st_level_31} \\    
    & \hspace{0.00cm} y^{sw,inv}_{l}, y^{line,inv}_{l} \in \{0,1\} ; \forall l \in \mathcal{L} \label{1st_level_32} \\
    & \hspace{0.00cm} y^{ha,inv}_{l,h} \in \{0,1\} ; \forall h \in \mathcal{H}_{l}, l \in \mathcal{L}^{ha}. \label{1st_level_33}
\end{align}

The objective function \eqref{1st_level_00} minimizes the annualized investment costs related to the installation of switching devices and new line segments, and their hardening. In addition, the objective function also accounts for operational costs over a selection of representative days $\mathcal{R}$, with $w_{r}$ indicating the frequency of $r \in \mathcal{R}$ within a year. The operational cost of each $r \in \mathcal{R}$ comprises expenses (switching, power imbalance) related to operating the system in a pre-contingency state during one selected period $t^{sp}$ (which can be defined, for instance, as the hour immediately before the start of the day or a peak-hour) combined with the worst-case expected cost of operating the system during the whole day while considering that failures might occur due to routine or wildfire events. Constraints \eqref{1st_level_01}--\eqref{1st_level_04} model the active and reactive nodal power balance for buses with substations and without substations. Constraints \eqref{1st_level_05}--\eqref{1st_level_06} model the voltage difference between adjacent buses as a function of the power flow through the line segment that connects them when this line segment is switched on, i.e., $z^{topo}_{l,r}=1$. Constraints \eqref{1st_level_07}--\eqref{1st_level_10} limit the active and reactive power flows following the linearization presented in \cite{Mashayekh2018}. Constraints \eqref{1st_level_11} impose bounds on nodal voltages, while the reference voltage for substations is set by constraints \eqref{1st_level_12}. Constraints \eqref{1st_level_13} and \eqref{1st_level_14} enforce power injection limits while \eqref{1st_level_15}--\eqref{1st_level_18} bound variables related to power imbalances. 

The group of constraints \eqref{1st_level_19}--\eqref{1st_level_33} models the behavior of the binary variables $z^{topo}$, $z^{sw}$, $y^{sw,inv}$, $y^{line,inv}$, and $y^{ha,inv}$. Variable $z^{topo}$ represents the line status, with $z^{topo} = 1$ if the line is on. If a line segment is switchable and turned off, then $z^{topo} = 0$. The same dynamics are considered for variables of candidate lines not created. Parameter $z^{init}$ represents the initial status of a line, with $z^{init} = 1$ if the line segment starts on or $z^{init} = 0$ otherwise. Variable $z^{sw}$ indicates turning a switchable line on/off, assuming $z^{sw} = 1$ if the action is made and $z^{sw} = 0$ otherwise. Variable $y^{sw,inv}_l$ represents investment decisions, with $y^{sw,inv}_l = 1$ prescribing the installation of a switching device on a given line segment $l \in \mathcal{L}^{sw,c}$, whereas variable $y^{line,inv}_l$ determines if a candidate line segment $l \in \mathcal{L}^{c}$ should be built or not. Finally, variable $y^{ha,inv}_{l,h}$ models investment decisions about making hardening investments $h \in \mathcal{H}_{l}$ on a given line segment $l \in \mathcal{L}^{ha}$, which might include different types of actions, such as line undergrounding, line coating, or even vegetation management. In summary, a given line can be part of one of six different categories: 1) existing and switchable; 2) existing and candidate to become switchable; 3) existing, non-switchable, and non-candidate to become switchable; 4) candidate that is switchable; 5) candidate that can become switchable if additional investment is made; and 6) candidate that cannot become switchable. Within this context, constraints \eqref{1st_level_19} ensure that line segments that are not switchable will not perform any switching action. Constraints \eqref{1st_level_20} express that line segments that are not a candidate to become switchable cannot receive investments to become switchable, while constraints \eqref{1st_level_21} ensure that the existing line segments should not be constructed again. Constraints \eqref{1st_level_22} enforce a permanent ``on'' status for non-switchable existing lines that are not candidates to receive switching devices. For each line segment candidate to become switchable, constraints \eqref{1st_level_23} determine that a switching action can only be performed if a switching device is installed. Analogously, for each line segment candidate to be built and to become switchable, constraints \eqref{1st_level_24} allow an investment in a switching device only if the solution indicates the construction of the line segment. Constraints \eqref{1st_level_25}--\eqref{1st_level_26} capture switching transitions for existing line segments. Constraints \eqref{1st_level_27} and \eqref{1st_level_28} model the behavior of status variable $z^{topo}_{l,r}$ for line segments that are candidates for construction. In this case, when a candidate line segment $l \in \mathcal{L}^{c} \setminus \mathcal{L}^{sw,e}$ is built without a switching device, $z^{topo}_{l,r} = 1, ~ \forall ~ r \in {\cal R}$. If the candidate line segment also has a switching device, variable $z^{topo}_{l,r}$ can be assigned either to 0 or 1. If line segment $l \in \mathcal{L}^{c}$ is not built, $z^{topo}_{l,r} = 0, ~ \forall ~ r \in {\cal R}$. Constraints \eqref{1st_level_29} avoid the simultaneous activation of line segments that can result in the formation of loops within the network, i.e., given a rule $k \in \mathcal{K}^{forb}$, at least one line segment of the set $\mathcal{L}^{forb}_{k}$ must be off ($z^{topo}=0$) to ensure a radial configuration for the distribution system. For any given grid, $\mathcal{L}^{forb}_{k}$ can be obtained via a depth-first search (DFS) algorithm \cite{Tarjan1972}. Constraints \eqref{1st_level_30} allow one hardening investment $h \in \mathcal{H}_{l}$ to be chosen for each line segment $l \in \mathcal{L}^{ha}$ that is a candidate for hardening. This type of investment directly affects the post-contingency operation by modifying the probability of a line failing given an adverse scenario (further explained in Section \ref{subsec:SecondLevel}). Constraints \eqref{1st_level_31}--\eqref{1st_level_33} express the binary nature of the status, switching and investment variables.

For didactic purposes, we hereinafter make use of the following concise representation of model \eqref{1st_level_00}--\eqref{1st_level_33}:
\begin{align}
    & \underset{{\substack{\boldsymbol{x}_{r}, \boldsymbol{z}_{r}, \boldsymbol{y}}}}{\text{Minimize}} 
    \hspace{0.40cm} \boldsymbol{c}^{inv^{\top}} \boldsymbol{y} + \sum_{r \in \mathcal{R}} w_{r} \bigg( \boldsymbol{c}^{sw^{\top}} \boldsymbol{z}_{r} + \boldsymbol{c}^{imb^{\top}} \boldsymbol{x}_{r} \notag \\
    & \hspace{3.00cm} + \sup_{\mathcal{Q} \in \mathcal{P}_{r}(\boldsymbol{f}_{r}^{p},\boldsymbol{y}^{ha})}\mathbb{E}_\mathcal{Q} \Big [ H_{r}(\boldsymbol{z}_{r},\boldsymbol{a}_{r}) \Big] \bigg) \label{1st_level_simplified_00} \\ \notag
    & \text{subject to:} \notag \\ 
    & \hspace{0.00cm} \boldsymbol{A} \boldsymbol{x}_{r} = \boldsymbol{b}_{r} ; \hspace{0.10cm} \forall r \in \mathcal{R} \label{1st_level_simplified_01} \\
    & \hspace{0.00cm} \boldsymbol{C} \boldsymbol{x}_{r} + \boldsymbol{L} \boldsymbol{z}_{r} \geq \boldsymbol{m}_{r} ; \hspace{0.10cm} \forall r \in \mathcal{R} \label{1st_level_simplified_02} \\
    & \hspace{0.00cm} \boldsymbol{N} \boldsymbol{z}_{r} + \boldsymbol{P} \boldsymbol{y} \geq \boldsymbol{u}_{r} ; \hspace{0.10cm} \forall r \in \mathcal{R} \label{1st_level_simplified_03} \\
    & \hspace{0.00cm} \boldsymbol{y}, \boldsymbol{z}_{r} \in \{0,1\}^{|\mathcal{L}|} ; \hspace{0.10cm} \forall r \in \mathcal{R}, \label{1st_level_simplified_04}
\end{align}

\noindent with $\boldsymbol{x}_{r} = \big[\boldsymbol{p}_{r}, \boldsymbol{q}_{r}, \boldsymbol{v}^{\dagger}_{r}, \boldsymbol{f}^{p}_{r}, \boldsymbol{f}^{q}_{r}, \boldsymbol{{\Delta}D}^{p+}_{r}, \boldsymbol{{\Delta}D}^{p-}_{r}, \boldsymbol{{\Delta}D}^{q+}_{r}, \boldsymbol{{\Delta}D}^{p-}_{r},$ $\boldsymbol{{\Delta}D}^{q-}_{r}\big]^{\top}$ and $\boldsymbol{z}_{r} = \big[\boldsymbol{z}^{topo}_{r}, \boldsymbol{z}^{sw}_{r}\big]^{\top}$ denoting respectively the variables for operative decisions and topology configuration of each representative day $r \in \mathcal{R}$, and $\boldsymbol{y} = \big[\boldsymbol{y}^{sw}, \boldsymbol{y}^{line}, \boldsymbol{y}^{ha}\big]^{\top}$ accounting for the (first-stage) investment decision variables. In the objective function \eqref{1st_level_simplified_00}, the first three terms stand for investments, switching, and imbalance costs, respectively. The last term represents the second-level problem (Subsection \ref{subsec:SecondLevel}) that maximizes the expected value of the post-contingency operation. Function $H_{r}\big(\boldsymbol{z}_{r}, \boldsymbol{a}_{r})$ denotes a third-level problem that minimizes the post-contingency operation for a given topology configuration ($\boldsymbol{z}_{r}$) and line availability scenario ($\boldsymbol{a}_{r}$) (explained in Subsection \ref{subsec:ThirdLevel}). Furthermore, the set of constraints \eqref{1st_level_simplified_01} represents the operative constraints \eqref{1st_level_01}--\eqref{1st_level_04} that models active and reactive power balance. The set of constraints \eqref{1st_level_simplified_02} represents the operative constraints \eqref{1st_level_05}--\eqref{1st_level_18} that models voltage differences, active and reactive power flow interaction, and operative variables limits. The set of constraints \eqref{1st_level_simplified_03} represents constraints \eqref{1st_level_19}--\eqref{1st_level_30} that models the behavior of binary variables for topology and investments decisions, and constraints \eqref{1st_level_simplified_04} represents constraints \eqref{1st_level_31}--\eqref{1st_level_33}, where the binary nature of variables $\boldsymbol{y}$ and $\boldsymbol{z}_{r}$ are defined.

\subsection{Multiperiod Operational Post-Contingency Problem} \label{subsec:ThirdLevel}

The critical task of planning/operating agents is to balance supply and demand while respecting the system's constraints by exploiting all available and scheduled resources. In the particular context of this work, to appropriately define the investments to be made in the distribution grid, this operational stage is represented as the third-level optimization problem under function $H_{r}$ for each representative day $r \in \mathcal{R}$. Formally, the model is similar to the first-level problem, except that a multiperiod operation is considered, and the topology information ($\boldsymbol{z}_{r}$) and line availability ($\boldsymbol{a}_{r}$) are known (inputs). The compact formulation for a given $r \in \mathcal{R}$ of this third-level problem is presented in \eqref{3rd_level_simplified_00}--\eqref{3rd_level_simplified_02}:
\begin{align}
    \hspace{-0.30cm} H_{r}(\boldsymbol{z}_{r},\boldsymbol{a}_{r}) = ~ & \underset{{\substack{\boldsymbol{x}_{t,r}^{c}}}}{\text{Minimize}} 
    \hspace{0.30cm} \frac{1}{|\mathcal{T}|} \sum_{t \in \mathcal{T}} \Big( \boldsymbol{c}^{oper^{\top}} \boldsymbol{x}_{t,r}^{c} \Big) \label{3rd_level_simplified_00} \\ \notag
    & \hspace{0.00cm} \text{subject to:} \notag \\ 
    & \hspace{0.00cm} \boldsymbol{A} \boldsymbol{x}_{t,r}^{c} = \boldsymbol{b}_{t,r} ; \hspace{0.10cm} \forall t \in \mathcal{T} \label{3rd_level_simplified_01} \\
    & \hspace{0.0cm} \boldsymbol{C} \boldsymbol{x}_{t,r}^{c} \geq \boldsymbol{m}_{t,r} - \boldsymbol{L} \boldsymbol{z}_{r} + \boldsymbol{V} \boldsymbol{a}_{r} \notag \\
    & \hspace{1.90cm} + \boldsymbol{W} (\boldsymbol{z}_{r} \pdot \boldsymbol{a}_{r}) ; \hspace{0.10cm} \forall t \in \mathcal{T}, \label{3rd_level_simplified_02}
\end{align}

\noindent where the set of variables $\boldsymbol{x}_{t,r}^{c}$ represents the variables for operative decisions of each period $t \in \mathcal{T}$ of a given representative day $r \in \mathcal{R}$. These operational decisions have the same role as the first-stage continuous decision variables in \eqref{1st_level_00}--\eqref{1st_level_33}, which includes power, voltage, line flow, and imbalance variables. Furthermore, as the operational stage \eqref{3rd_level_simplified_00}--\eqref{3rd_level_simplified_02} is aware of the contingency status of the line segments (i.e., $\boldsymbol{a}_{r}$ is an input), the model represents grid operation in a post-contingency state. Therefore, superscript $c$ is used to highlight that the variables in this stage are being optimized considering contingencies, thus differentiating them from the first-stage variables. 

In this problem, the objective function \eqref{3rd_level_simplified_00} minimizes the hourly average cost associated with electricity purchase and active and reactive power imbalance at every bus. Analogously to the first-stage constraints, \eqref{3rd_level_simplified_01} represents operative restrictions \eqref{1st_level_01}--\eqref{1st_level_04} of active and reactive power balance, for each period $t \in \mathcal{T}$; and the set of constraints \eqref{3rd_level_simplified_02} represents operative constraints \eqref{1st_level_05}--\eqref{1st_level_18} for each period $t \in \mathcal{T}$, modeling voltage differences, active and reactive power flow interaction, and operative variables limits. In this case, line availability will depend on the contingency variable $\boldsymbol{a}_{r}$ and on the first-stage topology configuration $\boldsymbol{z}_{r}$. Operator $\pdot$ stands for a Hadamard product, i.e., a component-wise multiplication between two vectors ($\boldsymbol{a}_{r}$ and $\boldsymbol{z}_{r}$).

\subsection{Decision-Dependent Ambiguity Set Modeling} \label{subsec:SecondLevel}

Following the discussion in the previous sections, the proposed investment planning methodology \eqref{1st_level_00}--\eqref{1st_level_33} considers the impact of operational decisions on the line availability probabilistic characterization with the following term in the objective function:
\begin{align}
    & \hspace{0.00cm} \sup_{\mathcal{Q} \in \mathcal{P}_{r}(\boldsymbol{f}_{r}^{p},\boldsymbol{y}^{ha})}\mathbb{E}_\mathcal{Q} \Big [ H_{r}(\boldsymbol{z}_{r},\boldsymbol{a}_{r}) \Big]. \label{2nd_level_00}
\end{align}

Problem \eqref{2nd_level_00} stands for the worst expected cost of the multiperiod post-contingency operation for a given representative day $r \in \mathcal{R}$ over the (ambiguity) set $\mathcal{P}_r(\boldsymbol{f}_{r}^{p},\boldsymbol{y}^{ha}) \in \mathcal{M}_{+}$, whose components are essentially different line availability probability distributions given the limited knowledge of failure probabilities, which are also dependent on the power-flow and on the hardening investment (first-stage) variables $(\boldsymbol{f}_{r}^{p}, \boldsymbol{y}^{ha})$. In the proposed methodology, the definition of an investment plan considers the usual nominal failure probability level associated with each line segment but also incorporates the fact that environmental conditions can create wildfire-prone circumstances. In this latter context, we argue that the power-flow levels $(\boldsymbol{f}_{r}^{p})$ impact the uncertainty characterization of the line availability. More specifically, the probability of a line starting a fire and becoming unavailable rises as the power through this line segment increases. Furthermore, as a counter-reaction, the hardening decisions $(\boldsymbol{y}^{ha})$ reduce line unavailability likelihood under extreme weather circumstances. Mathematically, the decision-dependent ambiguity set proposed in this work can be expressed as in \eqref{ambiguity_set_01}:
\begin{align}
    & \hspace{0.00cm} \mathcal{P}_{r}(\boldsymbol{f}_{r}^{p},\boldsymbol{y}^{ha}) = \Big\{ \mathcal{Q} \in \mathcal{M}_{+} (\mathcal{A}_{r}) ~ \Big| ~ \mathbb{E}_{\mathcal{Q}} \big[\boldsymbol{S} \hat{\boldsymbol{a}}_{r} \big] \leq \overline{\boldsymbol{\mu}}_{r}( \boldsymbol{f}_{r}^{p},\boldsymbol{y}^{ha}) \Big\}, \label{ambiguity_set_01}
\end{align}

\noindent where $\boldsymbol{S}$ is a matrix defined as $[\mathbb{I} ~ | ~ -\mathbb{I}]^{\top}_{2|\mathcal{L}| \times |\mathcal{L}|}$, with $\mathbb{I}$ being the identity matrix of size $|\cal L|$. The elements of vector $\hat{\boldsymbol{a}}_{r}= \mathbb{1} - \boldsymbol{a}_{r}$ are random variables associated with line unavailability, with $\mathbb{1}$ being a vector of ones and $\boldsymbol{a}_{r} \in \{0,1\}^{|\mathcal{L}|}$ indicating the availability of the line segments, i.e., each of its elements may be assigned to a value ${a}_{l,r} = 1$ if the corresponding line is available or ${a}_{l,r} = 0$ if it is unavailable. We assume that this contingency status stands for the whole operative period $\mathcal{T}$. Following standard industry and academic practices, the support of the random vector $\boldsymbol{a}_{r}$, is defined as:
\begin{align}
    & \hspace{0.00cm} \mathcal{A}_{r} = \bigg\{ \boldsymbol{a}_{r} \in \{0,1\}^{|\mathcal{L}|} ~ \bigg| ~ \sum_{l \in \mathcal{L}} a_{l,r} \geq |\mathcal{L}| - K \bigg\}, \label{ambiguity_set_02}
\end{align}

\noindent which essentially considers all failure events that involve at most $K$ line segments. Moreover, in the framework proposed in this work, the decision-dependent vector-function $\overline{\boldsymbol{\mu}}_{r}(\boldsymbol{f}_{r}^{p}, \boldsymbol{y}^{ha})$ is defined as:
\begin{align}
    & \hspace{-0.10cm} \overline{\boldsymbol{\mu}}_{r} ( \boldsymbol{f}_{r}^{p}, \boldsymbol{y}^{ha} ) = \boldsymbol{\gamma} + \text{diag} \Big( \boldsymbol{\beta}_{r} \pdot \boldsymbol{w}^{\beta} \big(\boldsymbol{y}^{ha}\big) \Big) \> \big|\boldsymbol{f}^{p}_{r}\big|, \label{ambiguity_set_03}
\end{align}

\noindent where vector $\boldsymbol{\gamma} = \left[ ~ \gamma_{l}^{\top} ~ | ~ \mathbb{0}^{\top} ~ \right]^{\top}_{1 \times 2|\mathcal{L}|}$ comprises an estimated nominal probability of failure associated with each line segment $l \in \cal{L}$ and $\mathbb{0}$ is a vector of zeros of size $|\mathcal{L}|$. This vector $\boldsymbol{\gamma}$ characterizes the exogenous nature of failure not related to any decision variable. Each component of $\boldsymbol{\beta}_r$ stands for the endogenous relationship between failure probability and power-flow levels through a specific line segment. Essentially, as the absolute value of the power flow through a line segment $l \in \cal{L}$ raises, the corresponding failure likelihood increases by a factor of ${\beta}_{l,r}$. In this context, hardening decisions taken in the planning stage decrease the impact of power-flow levels on failure probabilities. Formally, the effect of hardening decisions is expressed by the element-wise multiplication between $\boldsymbol{\beta}_r$ and $\boldsymbol{w}^{\beta}$ in \eqref{ambiguity_set_03}, with $\boldsymbol{w}^{\beta}$ defined as:
\begin{align}
    & \hspace{0.00cm} w^{\beta}_{l} \big(\boldsymbol{y}^{ha} \big) = 1 - \sum_{h \in \mathcal{H}_{l}} w^{ha}_{l,h} y^{ha,inv}_{l,h} ; \forall l \in \mathcal{L}^{ha} \label{ambiguity_set_04} \\
    & \hspace{0.00cm} w^{\beta}_{l} \big(\boldsymbol{y}^{ha} \big) = 1 ; \forall l \in \mathcal{L} \setminus \mathcal{L}^{ha}.\label{ambiguity_set_05}
\end{align}

% \noindent which takes into account the influence of the first-stage hardening decisions that can include different types of actions, such as line undergrounding, line coating, or even vegetation management (which is not technically a line hardening measure but can decrease failure probability by properly intervening in the surroundings of line segments). 

\section{Solution Methodology} \label{sec:Solution}

The two-stage formulation \eqref{1st_level_00}--\eqref{1st_level_33} proposed in Section \ref{sec:Problem} is intended to optimally determine grid investments combined with topology configuration to improve the flexibility of the system to deal with wildfire-prone conditions. In the previous section, we formulate this decision-making process as a trilevel optimization problem that cannot be directly solved via off-the-shelf solvers because of its complexity. Hence, in this section, we provide a tailored solution methodology to address this problem. To design this methodology, we leverage the convex properties of the innermost models to write an outer-approximation-based decomposition algorithm composed of tractable master- and sub-problems. In the next subsections, we thoroughly describe the procedures and reformulations needed to describe the proposed solution method.

\subsection{Master Problem} \label{subsec:Master}

Following the definition of the ambiguity set in \eqref{ambiguity_set_01}, for each representative day $r \in \mathcal{R}$, problem \eqref{2nd_level_00} can be expressed as the following mathematical programming formulation:
\begin{align}
    & \underset{{\mathcal{Q} \in {\mathcal{M}}^{+}}}{\text{Maximize}} 
    \hspace{0.20cm} \sum_{\boldsymbol{a}_{r} \in \mathcal{A}_{r}} H_{r}(\boldsymbol{z}_{r},\boldsymbol{a}_{r}) \mathcal{Q}(\boldsymbol{a}_{r}) \label{2nd_level_01} \\ \notag
    & \text{subject to:} \notag \\ 
    & \hspace{0.00cm} \sum_{\boldsymbol{a}_{r} \in \mathcal{A}_{r}} (S\hat{\boldsymbol{a}}_{r}) \mathcal{Q}(\boldsymbol{a}_{r}) \leq \overline{\boldsymbol{\mu}}_{r} ( \boldsymbol{f}_{r}^{p}, \boldsymbol{y}^{ha} ) :(\boldsymbol{\psi}_{r}) \label{2nd_level_02} \\
    & \hspace{0.00cm} \sum_{\boldsymbol{a}_{r} \in \mathcal{A}_{r}} \mathcal{Q}(\boldsymbol{a}_{r}) = 1 :(\varphi_{r}), \label{2nd_level_03}
\end{align}

\noindent which essentially identifies a probability measure $\mathcal{Q} \in {\mathcal{M}}^{+}$ that maximizes the expected cost of the post-contingency operation \eqref{2nd_level_00} while respecting \eqref{2nd_level_02} and \eqref{2nd_level_03}. Due to its convex properties, model \eqref{2nd_level_01}--\eqref{2nd_level_03} can be equivalently rewritten in its dual version as:
\begin{align}
    & \underset{{\psi_{l,r}, \varphi_{r}}}{\text{Minimize}} 
    \hspace{0.50cm} \boldsymbol{\psi}^{\top}_{r} \overline{\boldsymbol{\mu}}_{r} ( \boldsymbol{f}_{r}^{p}, \boldsymbol{y}^{ha} ) + \varphi_{r} \label{2nd_level_dual_01} \\ \notag
    & \text{subject to:} \notag \\ 
    & \hspace{0.00cm} \boldsymbol{\psi}^{\top}_{r} S \hat{\boldsymbol{a}}_{r} + \varphi_{r} \geq H_{r}(\boldsymbol{z}_{r}, \boldsymbol{a}_{r}) ; \forall \boldsymbol{a}_{r} \in \mathcal{A}_{r} \label{2nd_level_dual_02} \\
    & \hspace{0.00cm} \psi_{l,r} \geq 0 ; \forall l \in \big\{1,\dots,2 |\mathcal{L}|\big\} \label{2nd_level_dual_03} \\
    & \hspace{0.00cm} \varphi_{r} \in \mathbb{R}, \label{2nd_level_dual_04}
\end{align}

\noindent where variables $\psi$ and $\varphi$ are the dual variables associated with constraints \eqref{2nd_level_02} and \eqref{2nd_level_03}, respectively. Within this framework, problem \eqref{1st_level_00}--\eqref{1st_level_33} can be equivalently written as:
\begin{align}
    & \underset{{\substack{\boldsymbol{x}_{r}, \boldsymbol{z}_{r}, \boldsymbol{y}, \\
                           \boldsymbol{\psi}_{r}, \varphi_{r} }}}{\text{Minimize}} 
    \hspace{0.40cm} \boldsymbol{c}^{inv^{\top}} \boldsymbol{y} + \sum_{r \in \mathcal{R}} w_{r} \Big( \boldsymbol{c}^{sw^{\top}} \boldsymbol{z}_{r} + \boldsymbol{c}^{imb^{\top}} \boldsymbol{x}_{r} \notag \\
    & \hspace{4.00cm}  + \boldsymbol{\psi}^{\top}_{r} \overline{\boldsymbol{\mu}}_{r} ( \boldsymbol{f}_{r}^{p}, \boldsymbol{y}^{ha} ) + \varphi_{r} \Big) \label{1st_level_new_00} \\ \notag
    & \text{subject to:} \notag \\ 
    & \hspace{0.00cm} \text{Constraints \eqref{1st_level_simplified_01}--\eqref{1st_level_simplified_04}} \label{1st_level_new_01} \\
    & \hspace{0.00cm} \boldsymbol{\psi}^{\top}_{r} S \hat{\boldsymbol{a}}_{r} + \varphi_{r} \geq H_{r}(\boldsymbol{z}_{r}, \boldsymbol{a}_{r}) ; \hspace{0.10cm} \forall \boldsymbol{a}_{r} \in \mathcal{A}_{r}, r \in \mathcal{R} \label{1st_level_new_02} \\
    & \hspace{0.00cm} \psi_{l,r} \geq 0 ; \hspace{0.10cm} \forall l \in \big\{1,\dots,2 |\mathcal{L}|\big\}, r \in \mathcal{R} \label{1st_level_new_03} \\
    & \hspace{0.00cm} \varphi_{r} \in \mathbb{R} ; \hspace{0.10cm} \forall r \in \mathcal{R}, \label{1st_level_new_04}
\end{align}

\noindent Furthermore, constraints \eqref{1st_level_new_02} can be expressed equivalently as:
\begin{align}
    & \hspace{0.00cm} \varphi_{r} \geq \max_{\boldsymbol{a}_{r} \in \mathcal{A}_{r}}{ \Big\{ H_{r} (\boldsymbol{z}_{r}, \boldsymbol{a}_{r}) - \boldsymbol{\psi}^{\top}_{r} S \hat{\boldsymbol{a}}_{r} } \Big\} ; \forall r \in \mathcal{R}, \label{1st_level_new_05}
\end{align}

\noindent where the right-hand side can be relaxed and outer-approximated by cutting planes. This relaxed version of \eqref{1st_level_00}--\eqref{1st_level_33} is the master problem of the proposed decomposition algorithm whose approximation is iteratively improved by adding new cutting planes based on the solution of the sub-problems (Subsection \ref{subsec:Subproblem}). Note that the term $\boldsymbol{\psi}^{\top}_{r} \overline{\boldsymbol{\mu}}_{r} ( \boldsymbol{f}_{r}^{p}, \boldsymbol{y}^{ha} )$ in \eqref{1st_level_new_00} entails nonlinear products due to the product between elements of vectors $\boldsymbol{\psi}_{r}$ and $\boldsymbol{w}^{\beta} \big(\boldsymbol{y}^{ha}\big)$ with $|\boldsymbol{f}^{p}_{r}|$. We linearize these products via binary expansion and disjunctive constraints.

\subsection{Subproblem} \label{subsec:Subproblem}

The subproblem of the proposed decomposition algorithm is essentially the left-hand side of \eqref{1st_level_new_05}. Since $H_{r} (\boldsymbol{z}_{r}, \boldsymbol{a}_{r})$ is a convex minimization problem, we can equivalently replace it with its dual formulation and obtain a single-level maximization problem. By linearizing mixed-binary bilinear products, the resulting formulation is a Mixed-Integer Linear Programming (MILP) model that can be handled efficiently by MILP algorithms or off-the-shelf solvers. This subproblem provides the necessary information to build cutting planes for the previously described master problem.

\subsection{Outer-Approximation-Based Decomposition Algorithm} \label{subsec:Algorithm}

\begin{figure}[tb]
    \centering
    \includegraphics[width=0.45\textwidth,height=.5\textheight,keepaspectratio]{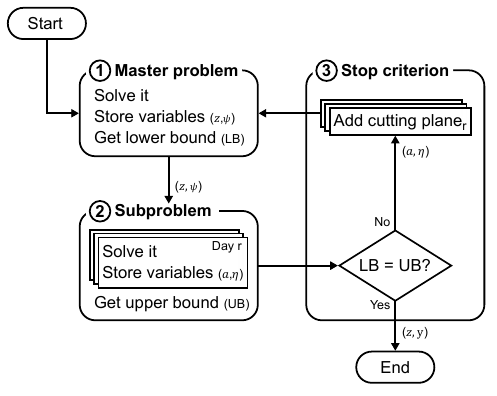}
    \caption{Solution algorithm.}
    \label{fig:solution_algorithm}
\end{figure}

The proposed decomposition algorithm, illustrated in Fig. \ref{fig:solution_algorithm}, iteratively solves the master and subproblems until the outer approximation of the right-hand side of \eqref{1st_level_new_05} properly describes it. At this point, the solution of the original investment planning problem \eqref{1st_level_00}--\eqref{1st_level_33} is obtained. To evaluate the closeness of the outer approximation, note that the optimal objective function value of the (relaxed) master problem provides a lower bound to the solution of \eqref{1st_level_00}--\eqref{1st_level_33}. Moreover, the solution of the subproblem can be combined with part of the objective function value of the master problem to compute an upper bound. Therefore, by comparing these upper- and lower-bounds, we can evaluate the accuracy of the outer approximation.

\section{Case Study} \label{sec:CaseStudy}

To illustrate the applicability of the proposed planning methodology, we use a 54-bus system based on \cite{Moreira2023}. We consider a 24-hour operation of 4 days representing a 1-year evaluation. The first two days represent operations during winter and spring, respectively, when we assume the wildfire risk is negligible. The third day corresponds to summer/fall also without wildfire risk, and the fourth day is associated with a summer/fall period when there is a high probability of wildfire occurrence. Furthermore, we assume $K = 1$ in \eqref{ambiguity_set_02} to characterize the support set $\mathcal{A}$. The proposed methodology is implemented in Julia$^{\text{\textregistered}}$ 1.6 on a server with one Intel$^{\text{\textregistered}}$ Core i7-10700K processor @ 3.80GHz and 64 GB of RAM, with Gurobi$^{\text{\textregistered}}$ 9.0.3. under JuMP$^{\text{\textregistered}}$.

\subsection{Experimental set up and input Data} \label{subsec:Input}

\begin{figure}[t]
\centerline{
    \subfigure[Representative days 1-3]{\includegraphics[width=0.24\textwidth]{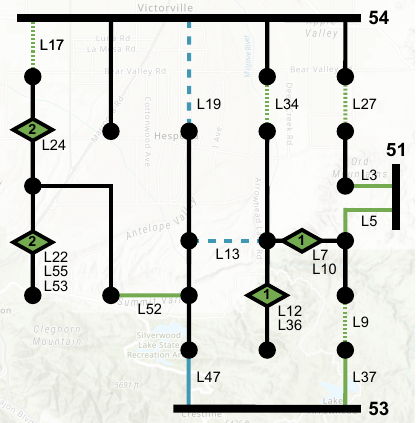}
    \label{fig:grid_123}}
    \subfigure[Representative day 4]{\includegraphics[width=0.24\textwidth]{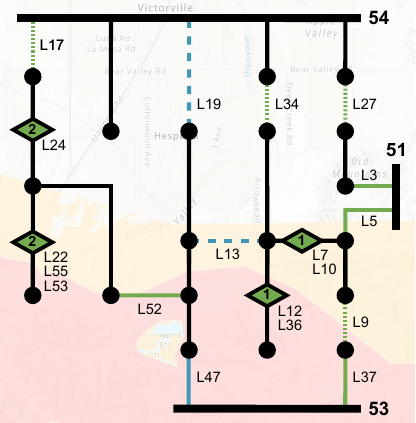}
    \label{fig:grid_4}}}
\caption{Simplified network representation of the 54-bus system with substations, selected buses, and lines characterized, respectively, by black rectangles, black circles, and line segments. The complete diagram can be accessed in \cite{DatasetPaper2024}.}
\label{fig:grid}
\end{figure}

We consider a distribution system with 54 buses and 57 lines. A simplified diagram of the test system is shown in Fig. \ref{fig:grid} for each of the 4 representative days under consideration, where we highlight substations 51, 53, and 54 (black rectangles), a few buses (black circles), and connections between them (line segments). The complete diagram can be accessed in \cite{DatasetPaper2024}, as well as all parameters used in this study. All lines that can either be switched or receive investment are labeled in the figure. Note, for example, that line segments L22, L55, and L53 are connected in sequence through intermediate buses and depicted in a simplified manner in Fig. \ref{fig:grid}. In addition, these lines are candidates to receive the hardening option 2. In this system, lines 19 and 13 are switchable lines that start open (blue dashed lines in the figure), while line 47 is a switchable line that starts closed (blue solid lines in the figure). Besides that, candidate lines 9, 17, 27, and 34 require an investment to be constructed (initially open, shown in the figure as green dotted lines). Existing lines 3, 5, 37, and 52 are candidates to become switchable (initially closed, shown in the figure as solid green lines). We also consider investment options for two different hardening actions for existing lines: (1) line undergrounding conversion for lines 7, 10, 12, and 36, and (2) line coating for lines 22, 24, 53, and 55. While the first hardening investment can reduce by 100\% the wildfire risk associated with power lines, the second option can reduce it by 60\% \cite{GridSafetyProgramSCE2018}. Therefore, in this experiment, we set to 1.0 and 0.6 the parameters that decrease the impact of power flow on failure probabilities based respectively on the investment options 1 and 2 ($w^{ha}$).

To illustrate the grid's wildfire-prone areas we rely on CPUC's High Fire-Threat District (HFTD) Map \cite{CPUC2024}. This map showcases the areas of California with an increased wildfire risk along with the associated utility company. For illustrative purposes, we consider the 54-system to be in San Bernardino County (Southern California). According to the HFTD map and following the system positioning, part of the grid is within the ``tier 2'' area, where an elevated risk is associated with overhead lines (depicted in yellow in Fig. \ref{fig:grid_4}). Another part of the grid would be inside the ``tier 3'' area, where there is an extreme risk associated with overhead lines (depicted in red in Fig. \ref{fig:grid_4}). Historically, fire season in Southern California lasts from late spring to early fall \cite{NWCG2024}, with its peak coinciding with the dry and warm summer season (June--August) \cite{Swain2021}. Recently, there has been observed a change in this pattern, with a longer duration and the peak moving sometimes to the fall season (September--November) \cite{Swain2021}. This change is caused by the rising temperatures and the shifting in precipitation seasonality, mostly due to climate change \cite{Swain2021}. Because of that, in our study, we assume that there is no wildfire risk associated with overhead lines during winter and spring, and therefore for representative days 1 and 2 we set $\beta$ equal to 0 (Fig. \ref{fig:grid_123}). For summer and fall, we consider a high level of wildfire risk during 50 days (representative day 4). In this case, we choose the $\beta$ values to represent a maximum failure probability of 90\% ($\beta \times \overline{F}$) for ``tier 3'' area and 60\% for ``tier 2'' area  (Fig. \ref{fig:grid_4}). We assume no wildfire risk for the rest of the days during summer and fall (representative day 3). For all the other lines outside the risk areas, we set $\beta$ equal to zero. Furthermore, for all line segments, we consider a nominal rate equal to 0.45 failure per year, which translates into a failure probability of 0.12\% for the next 24 hours, or 0.005\% for the next hour ($\gamma$).

Finally, in this study, we consider cost parameters mostly based on open data from the SCE utility company (in Southern California), using the company's unit cost guide for 2022 \cite{CostsSCE2022} and the report \cite{GridSafetyProgramSCE2018}. We also assume a switching cost of 100 \$ per hour, considering that the usage of switching devices is usually made manually by field crew \cite{Lei2018}. As for the energy cost, we assume the utility's residential rate of 01-Jan-24 \cite{SCE2024}. As the cost of loss of load depends on many factors, we consider a value of 2 \$/kWh according to the discussions on the CAISO's policy initiative \cite{Calpine2022}. The investment costs were updated for 2023 U.S. dollars and annualized considering a lifetime of 40 years. These annualized investment costs are shown in Table \ref{tab:costs}, as well as the costs for electricity and loss of load.

\begin{table}[!h]
    \centering
    \begin{threeparttable}
    \caption{Annualized investment, energy, and loss of load costs.}
    \renewcommand{\arraystretch}{1.30}
    \begin{tabular}{r|rl|c|}
        \textbf{} & \textbf{Cost} & & \textbf{Ref.} \\ \hline
        \textbf{New line connection}            &   1,618    & \$/km/year \tnote{a}         & \cite{CostsSCE2022} \\ 
        \textbf{Non-automated switching device} &     615    & \$/year per device \tnote{a} & \cite{CostsSCE2022} \\ 
        \textbf{Underground conversion}         & 232,026    & \$/km/year \tnote{a}         & \cite{GridSafetyProgramSCE2018} \\ 
        \textbf{Conductor covering}             &  33,102    & \$/km/year \tnote{a}         & \cite{GridSafetyProgramSCE2018} \\ 
        \textbf{Electricity}                    &       0.33 & \$/kWh                       & \cite{SCE2024} \\
        \textbf{Loss of load}                   &       2.00 & \$/kWh                       & \cite{Calpine2022} \\ 
    \end{tabular}
    \label{tab:costs}
    \begin{tablenotes}
        \item[a] The dollar year of the annual investment costs is 2023.
    \end{tablenotes}
    \end{threeparttable}
\end{table}

The demand for each day is hourly accounted for, with 24 periods. The demand value of each bus is constructed based on \cite{Moreira2023}, but the profile of the total demand was adapted considering the California ISO demand profile of days 04-Jan-2023, 05-Apr-2023, 05-Jul-2023, and 04-Oct-2023, for the representative days 1, 2, 3, and 4, respectively \cite{CAISO2024}. We choose the ``selected time period'' ($t^{sp}$) for the construction of the ambiguity set to be the hour where the operation has the highest total demand. The peak load on representative day 1 was at 7 pm, on day 2 at 9 pm, and on days 3 and 4 at 8 pm. Furthermore, we consider representative days 1 and 2 to account for 91 days each in a typical year, day 3 to represent 133 days, and day 4 to represent 50 days with fire risk. Therefore, the weight parameters $\boldsymbol{w}$ were set to $w_{1} = 2184$, $w_{2} = 2184$, $w_{3} = 3192$, and $w_{4} = 1200$, to represent the hours in a year for each representative day.

\subsection{Main Results} \label{subsec:Results}

To better analyze the proposed planning methodology, we first evaluate the model without the decision-dependent aspect, i.e., parameter $\beta$ set to zero (we, hereinafter refer to this case as ``without DDU''). Then, we run the model considering the configuration defined in Section \ref{subsec:Input} for the wildfire risk with $\beta$ parameter higher than zero (the ``with DDU'' case). We compare both results from the planning and switching decisions perspective. As shown in \cite{Moreira2023}, the ``with DDU'' case can leverage the cutting planes of the ``without DDU'' case as a warm-up phase. The key results of this comparison are presented in Table \ref{tab:main_results}.

\begin{figure}[t]
\centerline{
    \subfigure[Representative days 1-3]{\includegraphics[width=0.24\textwidth]{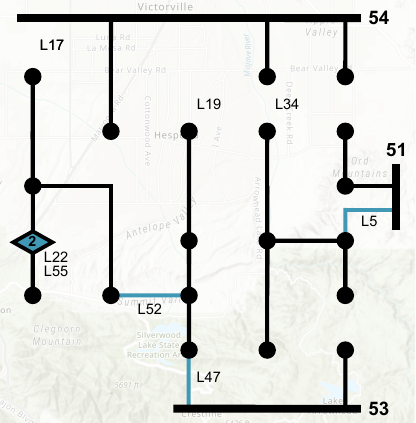}
    \label{fig:grid_results_123}}
    \subfigure[Representative day 4]{\includegraphics[width=0.24\textwidth]{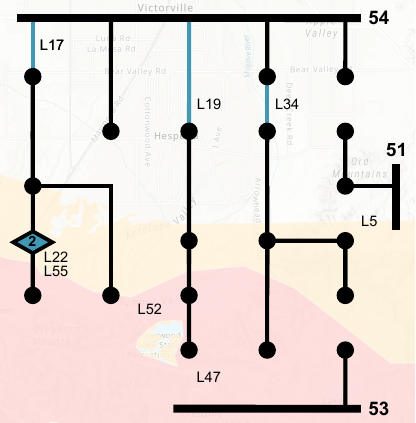}
    \label{fig:grid_results_4}}}
\caption{Final system topology for the ``with DDU'' case.}
\label{fig:grid_results}
\end{figure}

The ``with DDU'' case has as a solution a portfolio of different investments, prescribing the construction of lines 17 and 34, and turning lines 5 and 52 to be switchable. Hardening investments were made only in the ``with DDU'' case, where lines 22 and 55 with the hardening investment option 2 (line coating). These investments allow for topology reconfiguration for representative day 4 when there is wildfire risk. The final topology for days 1-3 and day 4 can be seen in Fig. \ref{fig:grid_results} for the ``with DDU'' case. As the ``without DDU'' case did not perform any topology reconfiguration, the final grid is equal to Fig. \ref{fig:grid_results_123}. It can be noted that, if no switching actions are performed, the buses inside the wildfire-prone area are entirely supplied by substations 51 and 53. When the topology reconfiguration is performed (Fig. \ref{fig:grid_results_4}), these buses have their energy mostly supplied by bus 54. In this configuration, these buses become the end of their branches, with a reduced power flow through the lines inside this high-threat area, and, consequently, a decrease in the wildfire risk. In Figure \ref{fig:grid_results}, it can also be seen that lines 22 and 55 received investment to be hardened. As they are at the end of their branch, no topology reconfiguration would allow a reduction of power flow levels, therefore, it is worth investing in hardening.

\begin{table}[!ht]
    \centering
    \begin{threeparttable}
    \footnotesize
    \newcolumntype{R}[1]{>{\raggedleft\arraybackslash}p{#1}}
    \caption{Main results}
    \renewcommand{\arraystretch}{1.30}
    \begin{tabular}{rc|R{2.0cm}|R{2.4cm}|}
        \multicolumn{1}{l}{\textbf{}} & \multicolumn{1}{l|}{\textbf{}} & \multicolumn{1}{c|}{\textbf{W/out DDU}} & \multicolumn{1}{c|}{\textbf{With DDU}} \\ \hline
        \multicolumn{2}{r|}{\textbf{New lines}} & - & 17; 34 \\ \hline
        \multicolumn{2}{r|}{\textbf{New switchable lines}} & - & 5; 52 \\ \hline
        \multicolumn{1}{r|}{\multirow{4}{*}{\makecell[r]{\textbf{Sw. actions} \\ \textbf{(line index)}}}} & \textbf{Day 1} & - & - \\
        \multicolumn{1}{r|}{} & \textbf{Day 2} & - & - \\
        \multicolumn{1}{r|}{} & \textbf{Day 3} & - & - \\ 
        \multicolumn{1}{r|}{} & \textbf{Day 4} & - & 5; 17; 19; 34; 47; 52 \\ \hline
        \multicolumn{1}{r|}{\multirow{2}{*}{\makecell[r]{\textbf{Hardening} \\ \textbf{investments}}}} & \textbf{1} & - & - \\
        \multicolumn{1}{r|}{} & \textbf{2} & - & 22; 55 \\ \hline
        \multicolumn{2}{r|}{\textbf{Annual investment costs (\$)}} & 0 & 91,049 \\ \hline
        \multicolumn{2}{r|}{\textbf{Time (min)}} & 27.5 & 65.8 \tnote{a} 
    \end{tabular}
    \label{tab:main_results}
    \begin{tablenotes}
      \item[a] The model ``with DDU'' needed 38.3 min to run. As it considered the ``without DDU'' cutting planes that lasted 27.5 min, in total the model needed 65.8 min.
    \end{tablenotes}
    \end{threeparttable}
\end{table}

We also analyze both cases from a financial viewpoint. The solution for the ``with DDU'' case would require an investment of \$ 1,230, \$ 3,753, \$ 86,066, respectively, for switching devices, line construction, and line hardening investment, summing up \$ 91,049 per year. On the other hand, the ``without DDU'' case has no investment or switching action prescribed, as there is no incentive to change the topology.

\subsection{Out-of-sample Analysis} \label{subsec:Outofsample}

With the planning and switching decisions (grid topology and investments) defined for both ``with DDU'' and ``without DDU'' cases, we perform an out-of-sample analysis using Monte Carlo simulation \cite{Mooney1997}. In the simulation, a one-year operation is considered as one scenario and we simulate 500 scenarios. By doing this analysis, we can verify the value and impact of considering the DDU aspect for the distribution system planning and operation. Formally, for each year/scenario, the following 4 steps are taken. \textbf{Step 1}: We calculate the power flow of each line in the operation of the distribution grid without any contingencies for each day in the year. \textbf{Step 2}: We calculate the probability of failure of each line for every hour and day with expression \eqref{ambiguity_set_03}. For both cases, ``with DDU'' and ``without DDU'', we calculate this failure with the $\beta$ used in the ``with DDU'' case. Although in the construction of the ambiguity set, we only consider the failure associated with the power flow at the selected period $t^{sp}$, in this out-of-sample analysis, we calculate individually the failure rate of each period $t \in \mathcal{T}$. \textbf{Step 3}: We randomly sample the status of each line in each hour and day based on the calculated probability of failure using a Bernoulli distribution, with a ``success'' outcome indicating out-of-service. \textbf{Step 4}: We evaluate the impact of the contingencies in the operation of each case, in each day. First, we calculate the total active power deficit proportion compared to the entire demand of the whole year. Then, we calculate two reliability metrics, SAIDI (System Average Interruption Duration Index) and SAIFI (System Average Interruption Frequency Index). These two metrics are standardized by IEEE specifically for distribution systems \cite{IEEE2022}. %While the first ``indicates the total duration of interruption for the average customer during a predefined period'', the second ``indicates how often the average customer experiences a sustained interruption over a predefined period''. 

The results of all scenarios show an average loss of active load of 5.5\% of the total demand in a year for the ``without DDU'' case and a value of 0.4\% for the ``with DDU'' case. This difference can also be perceived when we evaluate the loss of load in terms of CVaR$_{95\%}$. While the ``without DDU'' case has a CVaR$_{95\%}$ of 5.6\%, the ``with DDU'' case has a CVaR$_{95\%}$ of 0.4\%. The loss of load observed in both cases can be summarized in monetary values considering a value of loss of load equal to 2.00 \$/kWh. On average, the cost of the deficit for a year would be \$ 5,804,245 and \$ 384,486 when considering the ``without DDU'' and ``with DDU'' cases, respectively. Furthermore, as previously discussed, the operation of the ``with DDU'' case would require an annual investment of \$ 91,049, and the ``without DDU'' case would not require any investment. This cost of investment pays out since the difference in the average expected cost of loss of load between cases is \$ 5,419,759. 

Finally, regarding the reliability metrics, the customers of this grid would expect, on average, 280 hours with no energy or 11.7 days in a year (SAIDI) if the grid was planned following the ``without DDU'' case. In this context, these customers would expect, on average, 59 interruptions throughout the year (SAIFI). On the other hand, if the grid was planned following the decisions indicated by the ``with DDU'' case, customers would expect, on average, 18 hours (SAIDI) with no energy and 14 interruptions (SAIFI) throughout the year. In addition, while the ``without DDU'' case has a CVaR$_{95\%}$ of SAIDI equal to 284 hours and a CVaR$_{95\%}$ of SAIFI equal to 62 interruptions, the ``with DDU'' case has a CVaR$_{95\%}$ of SAIDI equal to 19 hours and a CVaR$_{95\%}$ of SAIFI equal to 16 interruptions. For comparison purposes, we also evaluated the model with a counterfactual of 500 scenarios assuming no wildfire risk for the representative day 4. On average the results for the deficit, SAIDI, and SAIFI are, respectively, 0.02 \% of total demand, 1.2 hours, and 0.5 interruptions.
 
\subsection{Sensitivity to the Weight of Representative Days} \label{subsec:Sensitivity}

We conduct a sensitivity analysis to assess the impact of the weights assigned to representative day 4 on the decisions made by the proposed methodology. We consider five contexts: 10, 20, 30, 40, and 50 days of wildfire for representative day 4. We use the same weight for representative days 1 and 2 as before and adjust representative day 3 accordingly. Table \ref{tab:main_results_sensitivity} summarizes the main results for the ``with DDU'' case. As the ``without DDU'' case still does not indicate any investment regardless of the weight of representative day 4, its results are not included in Table \ref{tab:main_results_sensitivity}.

\begin{table}[!ht]
    \centering
    \begin{threeparttable}
    \footnotesize
    \newcolumntype{R}[1]{>{\raggedleft\arraybackslash}p{#1}}
    \caption{Sensitivity analysis results}
    \renewcommand{\arraystretch}{1.30}
    \begin{tabular}{rc|R{0.73cm}|R{0.73cm}|R{0.73cm}|R{0.73cm}|R{0.73cm}|}
        \multicolumn{2}{r|}{} & \multicolumn{5}{c|}{\textbf{N. of days with wildfire risk in a year}} \\ \cline{3-7}
        \multicolumn{1}{l}{\textbf{}} & \multicolumn{1}{l|}{\textbf{}}& \multicolumn{1}{c|}{\textbf{10}} & \multicolumn{1}{c|}{\textbf{20}} & \multicolumn{1}{c|}{\textbf{30}} & \multicolumn{1}{c|}{\textbf{40}} & \multicolumn{1}{c|}{\textbf{50}} \\ \hline
        \multicolumn{2}{r|}{\textbf{New lines}} & 17; 34 & 17; 34 & 17; 34 & 17; 34 & 17; 34 \\ \hline
        \multicolumn{2}{r|}{\textbf{New switchable lines}} & 5; 52 & 5; 52 & 5; 52 & 5; 52 & 5; 52 \\ \hline
        \multicolumn{1}{r|}{\multirow{2}{*}{\makecell[r]{\textbf{Hardening} \\ \textbf{investments}}}} & \textbf{1} & - & - & - & - & - \\
        \multicolumn{1}{r|}{} & \textbf{2} & - & 22 & 22 & 22 & 22; 55 \\ \hline
        \multicolumn{2}{r|}{\textbf{Annual invest. costs (\$)}} & 4,983 & 48,347 & 48,347 & 48,347 & 91,049 \\ \hline
        \multicolumn{2}{r|}{\textbf{Time (min)}} & 49.9 & 50.0 & 57.3 & 57.7 & 66.5
    \end{tabular}
    \label{tab:main_results_sensitivity}
    \end{threeparttable}
\end{table}

Following Table \ref{tab:main_results_sensitivity} it can be seen that, as the number of days with wildfire risk increases, the model decides to make more investments, as expected. In all cases, the corresponding solutions comprise investments in upgrading the grid by building line segments 17 and 34 and by turning lines 5 and 52 to be switchable. However, as the number of days with wildfire risk increases, hardening investment starts to take place and increases, while (i) with 10 days, no hardening option is chosen; (ii) with 30 days line 22 is upgraded with hardening option 2; (iii) and with 50 days, lines 22 and 55 are enhanced with this same hardening option. It can also be noted that the computing time to run the methodology also tends to increase as the wildfire risk increases.

\section{Conclusion} \label{sec:Conclusion}

We developed in this paper a DDU-aware model for planning distribution systems considering the risk of wildfire and the relationship between power line operation and fire ignition. The investment decisions (made in the first stage) of installing new lines and switching devices can increase the flexibility of the grid operation to change the topology on wildfire-prone days (in the second stage). By taking these actions, the power flow levels through high-threat areas can be decreased, therefore also reducing the lines' failure probability, which is modeled as dependent on the power flow itself and on the environmental conditions. The model also considers the investment option of hardening lines, where the user can input different hardening options with different costs and impacts on the considered relationship between power line operation and wildfire ignition. Since the model uses a DRO framework and has three levels, we also designed an iterative methodology to solve it.

Through a 54-bus case study, we evaluated the effectiveness of the proposed methodology. The results showed that, considering or not the DDU aspect leads to different investment decisions. When the DDU aspect is ignored, the wildfire impacts are underestimated and no investment decision and topology change is made. When the DDU aspect is considered, investments are proposed, and the grid is more prepared to deal with the wildfire risk. An out-of-sample analysis attested that when no investments are made, a higher cost of lost load is expected. In fact, the cost of loss of load associated with neglecting the DDU aspect is also much higher than the amount of funds needed to implement the investment plan obtained while taking DDU into account. Besides that, the out-of-sample analysis showed that on average the SAIFI and SAIDI for this numerical experiment are, respectively, 59 interruptions and 280 hours for the model ``without DDU'', and 14 interruptions and 18 hours for the model ``with DDU''.

%This line of research can be enhanced or adapted in different aspects. For future work, we suggest considering the operation of batteries and DERs for microgrid formation, modeling the switching actions as recourse decisions and correlated cascading failures, and considering other ways to linearize the problem other than binary expansion.

% \begin{enumerate}
%     \item Adapt the framework to different topological conditions such as distribution systems operated as meshed networks and transmission systems;
%     \item Consider the operation of batteries and DERs for microgrid formation;
%     \item Model the switching actions as recourse decisions and correlated cascading failures;
%     \item Consider other ways to linearize the problem. The binary expansion approach increased significantly the number of binary variables needed, which increased the computational burden;
%     \item Investigate other ways to construct the ambiguity set. We modeled the line failure probability of the whole multiperiod operation based on a single period, which can be unsuitable for some cases.
% \end{enumerate}

% === Bibliography === %
\bibliographystyle{IEEEtran}
\bibliography{References}

\end{document}